\documentclass[12pt]{article}

\usepackage{amssymb,amsmath,amsthm}                 %amsthm adds "." after theorem numbers.
\usepackage{bm}                                     %+ bold math

\def\liml{\mathop{\lim}\limits}                     %+
\def\suml{\mathop{\sum}\limits}                     %+
\def\eq#1{\begin{equation}#1\end{equation}}         %+
\def\eqs*#1{\begin{eqnarray*}#1\end{eqnarray*}}     %+
\def\eqss#1{\begin{eqnarray}#1\end{eqnarray}}       %+
\DeclareMathSymbol{\ell}{\mathord}{letters}{96}     %+ to preserve \ell for arXiv
\def\cdc{,\ldots,}                                  %+
\def\1n{1\cdc n}                                    %+
\newcommand{\R}{\mathbb{R}}                         %+
\def\a{\alpha}                                      %+
\def\g{\gamma}                                      %+
\def\ovec{\overrightarrow}                          %+
\def\FF{\mathop{\cal F}\nolimits}                   %+
\def\Fij{\FF_{ij}}                                  %+
\def\Fjj{\FF_{jj}}                                  %+
\def\Fii{\FF_{ii}}                                  %+
\def\ainv{\tfrac{1}{\a}}                            %+
\def\cd{\!\cdot}                                    %+
\def\_#1{\mathop{\hspace{-0.41ex}^{}_{#1}}}         %+
\def\proof{{\noindent\bf Proof. }}                  %+
\def\atn{\a^{\frac{2}{n}}}                          %+
\def\squareforqed{\hbox{\rlap{$\sqcap$}$\sqcup$}}   %+
\def\qed{\ifmmode\squareforqed\else{\unskip\nobreak\hfil\penalty50\hskip1em\null\nobreak\hfil\squareforqed  %...
\parfillskip=0pt\finalhyphendemerits=0\endgraf}\fi} %+
\newtheorem{thm}{Theorem}{\bfseries}{\itshape}      %+
\newtheorem{prop}{Proposition}{\bfseries}{\itshape} %+
\newtheorem{corol}{Corollary}{\bfseries}{\itshape}  %+
\newtheorem{defin}{Definition}{\bfseries}{\upshape} %+
\def\SS{{\cal S}}

\def\Ff{\mathsf F}                                  %forest
\def\Ffb{\mathsf{\bar F}}                           %forest
\def\e{\mathrm e}                                   %+
\def\ssi{\hspace{.1em}}
\def\Up#1{\vspace{-#1em}}
\providecommand{\url}[1]{#1}
\csname url@samestyle\endcsname

\title{A Class of Graph-Geodetic Distances Generalizing %both
       \\the Shortest-Path and the Resistance Distances
}

\author{Pavel Chebotarev\footnote{E-mail: {\tt chv@member.ams.org}}\\%, pavel4e@gmail.com}}\\
{\normalsize Institute of Control Sciences of the Russian Academy of Sciences}\\
{\normalsize 65 Profsoyuznaya Street, Moscow 117997, Russia}%\\
}
\date{}

\begin{document}
\maketitle
%\unitlength 1.50mm

\begin{abstract}
A new class of distances for graph vertices is proposed.
This class contains parametric families of distances which reduce to the shortest-path, weighted shortest-path, and the resistance distances at the limiting values of the family parameters. The main property of the class is that all distances it comprises are graph-geodetic: $d(i,j)+d(j,k)=d(i,k)$ if and only if every path from $i$ to $k$ passes through~$j$.
The construction of the class is based on the matrix forest theorem and the transition inequality. %for spanning rooted forests.

\bigskip

%15A51%, %Stochastic matrices
% 15A09, %Matrix inversion, generalized inverses,
% 60J10, %Markov chains with discrete parameter
% 60J22, %Computational methods in Markov chains
\end{abstract}

\section{Introduction}
\label{sec_intro}

The classical distance for graph vertices is the shortest path distance~\cite{BuckleyHarary90}. Another distance, which is almost classical, is the resistance distance~\cite{GvishianiGurvich87,KleinRandic93,Gurvich10}, which is proportional to the commute-time distance \cite{Nash-Williams59,GobelJagers74,ChandraRaghavan89}. 

The forest distances $\tilde d_\a(i,j)$ \cite{CheSha01} form a one-parametric family %of graph distances
converging to the discrete distance as $\a\to0$ ($\tilde d_0(i,j)=1$ whenever vertices $i$ and $j$ are distinct) and becoming proportional to the resistance distance as $\a\to\infty$. The parameter $\a$ controls the relative influence of short and long paths connecting two vertices on the distance between them.

In a recent paper \cite{YenSaerensShimbo08} (see also \cite{MantrachYenCallut09}), the authors construct a parametric family of graph dissimilarity measures whose extrema are the weighted shortest path distance and the resistance distance. It is noteworthy that in clustering tasks, the best performance is obtained with intermediate values of the family parameter. On the other hand, the corresponding intermediate measures
break the triangle inequality, so they need not be distances (in this paper, we use the term ``distance'' in the sense of a metric space).

Thus, there is a demand in certain applications (these include data analysis, computer science, mathematical chemistry and some others) for a class of graph \emph{distances\/} whose extreme properties are similar to those of the dissimilarity measures proposed in~\cite{YenSaerensShimbo08}. Such a class is introduced in this paper. It comprises logarithmically transformed forest distances, and its construction is based on the matrix forest theorem \cite{CheSha97} and the transition inequality~\cite{Che10bottl}. The logarithmic transformation not only leads
to the   shortest-path /
weighted shortest-path distance at $\a\to 0^+$ and
to the     resistance  distance at $\a\to \infty$, but also, for every $\a>0$, it ensures the remarkable \emph{graph-geodetic property}: $d(i,j)+d(j,k)=d(i,k)$ if and only if every path from $i$ to $k$ passes through~$j$.

We now introduce the necessary notation.
Let $G$ be a weighted multigraph (a weighted graph, where multiple edges are allowed) with vertex set $V(G)=\{\1n\}$, $n>1$ and edge set~$E(G)$. We assume that $G$ has no loops. For $i,j\in V(G)$, let $n_{ij}\in\{0,1,\ldots\}$ be the number of edges incident to both $i$ and $j$ in~$G$; for every $p\in\{\1n_{ij}\}$, $w_{ij}^p>0$ is the weight of the $p\/$th edge of this type; let $w_{ij}=\sum_{p=1}^{n_{ij}}w_{ij}^p$ (if $n_{ij}=0$, we set $w_{ij}=0$) and $W=(w_{ij})_{n\times n}$. $W$ is the symmetric \emph{matrix of total edge weights} of~$G$. %for all pairs of vertices.

\smallskip
A \emph{rooted tree\/} is a connected and acyclic weighted graph in which one vertex, called the {\it root}, is marked. A~{\it rooted forest\/} is a graph, all of whose connected components are rooted trees. The roots of those trees are, by definition, the roots of the rooted forest.

\smallskip
By the weight of a weighted graph $H$, $w(H)$, we mean the product of the weights of all its edges. If $H$ has no edges, then $w(H)=1$. The weight of a set $\SS$ of graphs, $w(\SS)$, is the total weight of the graphs belonging to~$\SS$; the weight of the empty set is zero. If the weights of all edges are unity, i.e.\ the graphs in $\SS$ are actually unweighted, then $w(\SS)$ reduces to the cardinality of~$\SS$.

\smallskip
For a given weighted multigraph $G$, by $\FF=\FF(G)$, $\Fij=\Fij(G)$, and $\Fij^{(p)}=\Fij^{(p)}(G)$ we denote the set of all spanning rooted forests of $G$, the set of all forests in $\FF$ that have vertex $i$ belonging to a tree rooted at~$j$, and the set of all forests in $\Fij$ that have exactly $p$ edges. Let
\vspace{-.8em}
\eq{
\label{e_fij}%+
f=w(\FF),\quad f_{ij}=w(\Fij),\;\; \text{and}\;\;\, f_{ij}^{(p)}=w\big(\Fij^{(p)}\big),\quad i,j\in V(G),\quad 0\le p<n; }
by $F$ we denote the matrix $(f_{ij})_{n\times n}$; $F$~is called the \emph{matrix of forests\/} of~$G$.

\smallskip
Let $L=(\ell_{ij})$ be the Laplacian matrix of $G$, i.e.,
\eqs*{
%\label{e_lij}
\ell_{ij}=
        \begin{cases}
        -w_{ij},               &j\ne i,\\
         \suml_{k\ne i}w_{ik}, &j  = i.
        \end{cases}
}

Consider the matrix
$$
Q=(q_{ij})=(I+L)^{-1}.
$$
By the matrix forest theorem\footnote{Cf.\ Theorems 1 to 3 in \cite{ChungZhao10}.} \cite{CheSha95,CheSha97,CheAga02ap}, for any weighted multigraph $G$, $Q$ does exist and
\eq{
\label{e_mft}%+
q_{ij}=\frac{f_{ij}}{f},\quad i,j=\1n.
}
Consequently, $F=fQ=f\cd\!(I+L)^{-1}$ holds.
$Q$~can be considered as a matrix providing a proximity (similarity) measure for the vertices of~$G$\/ \cite{CheSha97,Che08DAM}.

\smallskip
By $d^s(i,j)$ we denote the \emph{shortest path distance},\footnote{The weighted shortest path distance will be considered in Section~\ref{s_weighted}.} i.e., the number of edges in a shortest path between $i$ and $j$ in~$G$; by $d^r(i,j)$ we denote the \emph{resistance distance\/} between $i$ and $j$ defined as follows:
\eq{
\label{def_resdist}%+
d^r(i,j)=\ell_{ii}^++\ell_{jj}^+-2\ell_{ij}^+,
}
where $(\ell_{ij}^+)_{n\times n}=L^+$ is the Moore-Penrose generalized inverse of the Laplacian matrix $L$ of~$G$.
$d^r(i,j)$ is equal to the effective resistance between $i$ and $j$ in the resistive network whose line conductances equal the edge weights $w_{ij}^p$ in~$G$.
If $G$ is connected, then\footnote{In fact, for a connected graph, $L^+=(L+\a\bar J)^{-1}-\a^{-1}\bar J$ with any $\a\ne0$ (Propositions~7 and~8 in \cite{CheSha98}, where the more general case of a multicomponent graph is considered). This expression with $\a=n$ is presented in \cite[page~88]{KleinRandic93}. For other related references, see Remarks on Proposition~15 in~\cite{CheAga02ap}.}
\eq{
\label{e_L+J}%+
L^+=(L+\bar J)^{-1}-\bar J,
}
where $\bar J$ is the $n\!\times\!n$ matrix with all entries $\frac{1}{n}$. Furthermore, by \cite[Theorem~3]{CheSha98}
$$
\ell_{ij}^+=\frac{f_{ij}^{(n-2)}-\tfrac{1}{n}f^{(n-2)}}{nt}, \quad i,j\in V(G)
$$
holds, where $f^{(n-2)}$ is the total weight of spanning rooted forests with $n-2$ edges and $t$ is the total weight of spanning trees in~$G$. By virtue of \eqref{def_resdist} this yields
\begin{corol}[{to Theorem~3 of \cite{CheSha98} and \eqref{e_L+J}}]
\label{cor_fores}
If\/ $G$ is connected$,$ then
\eq{
\label{eq_cor_fores}%+
d^r(i,j)=x_{ii}+x_{jj}-2x_{ij}=\frac{f_{ii}^{(n-2)}+f_{jj}^{(n-2)}-2f_{ij}^{(n-2)}}{nt},\quad i,j\in V(G),
}
where $(x_{ij})=(L+\bar J)^{-1}$.
\end{corol}

%If $G$ is unweighted, then
The forest representation \eqref{eq_cor_fores} is a counterpart of the classical 2-tree expression for $d^r(i,j)$ (see, e.g., \cite[Theorem~7-4]{SeshuReed61} and \cite{Shapiro87}); %,Chan69
it will be of use in Section~\ref{s_asymp}.
%-%result obtained in Section~5 of~\cite{Bapat99RD}.
%Corollary~\ref{cor_fores} will be of use in Section~\ref{s_asymp}.

\medskip
In Section~\ref{s_class} we introduce a new class of intrinsic graph distances and in Sections~\ref{s_bottle}--\ref{s_weighted} we study its properties.

\section{Logarithmic forest distances}
\label{s_class}

Suppose that $G$ is a connected weighted multigraph. Let
\vspace{-.3em}
\eq{
\label{Qa}%+
Q_\a=(I+L_\a)^{-1},
}
where $\a$ is a real parameter, $I$ is the identity matrix, and $L_\a$ is the Laplacian matrix of the multigraph $G_\a$ resulting from $G$ by a certain transformation of edge weights. This transformation generally depends on $\a$; for example, if every edge weight is multiplied by ${\a>0,}$ then\footnote{In this case, \eqref{Qa} is called the \emph{regularized Laplacian kernel\/} of $G$ with diffusion factor $\a$ (see \cite{CheSha97,SmolaKondor03,ShimboItoM09}).} $L_\a=\a L$, where $L$ is the Laplacian matrix of~$G$.

Define the matrix $H_\a$ as follows:
\vspace{-.5em}
\eq{
\label{Ha}%+
H_\a=\g\,(\a-1)\,\ovec{\log_\a Q_\a},
}
where $\a\!>\!0,$ $\a\!\ne\!1$, $\g$ is a positive factor, %independent of $G$,
and $\ovec{\varphi(Q_\a)}$ with $\varphi$ being a function stands for elem\-entwise operations, i.e., operations applied to each entry of $Q_\a$ separately.
Finally, consider
\vspace{-.2em}
\eq{
\label{Da}%+
D_\a=\tfrac{1}{2}(h_\a{\bm1}'+\bm1 h'_\a)-H_\a,
}
where $h_\a$ is the column vector containing the diagonal entries of $H_\a$, $h'_\a$ is the transpose of $h_\a$, ${\bm1}$ and ${\bm1}'$ being the column of $n$ ones and its transpose. The elementwise form of \eqref{Da} is: $d_{ij}(\a)=\frac{1}{2}(h_{ii}(\a)+h_{jj}(\a))-h_{ij}(\a),\,$ $i,j=\1n.$ This is a standard transformation used to obtain a metric from a symmetric similarity measure (see, e.g., the inverse covariance mapping in~\cite{DezaLaurent97}). As Theorem~\ref{th_dist} below states, $D_\a$ determines a metric on %is a matrix of distances between
the vertices of~$G$.

\smallskip
Since $\liml_{\a\to1}\left((\a-1)/\ln\a\right)=1$, we extend Eq.\,\eqref{Ha} to $\a=1$ as follows:
\vspace{-.5em}
\eq{
\label{Ha1}%+
H_{1}=\g\,\ovec{\ln Q_1},
}
which preserves continuity.
This extension %used for the definition of $D_1$ according to \eqref{Da}
is assumed throughout the paper.

\begin{thm}
\label{th_dist}
For any connected multigraph $G$ and any $\a,\g>0,$ $D_\a=(d_{ij}(\a))_{n\times n}$ defined by Eqs.\ \eqref{Qa}--\eqref{Ha1} is a matrix of distances on\/~$V(G)$.
\end{thm}

Before proving Theorem~\ref{th_dist} we represent the entries of $D_\a$ in terms of the weights of spanning forests in~$G_\a$.
Let
\eq{
\label{e_fija}%+
f_{ij}(\a)=w(\Fij(G_\a)),\quad i,j=\1n.
}

\smallskip
\begin{prop}
\label{p_d1}
For any connected multigraph $G$ and any $\a,\g\!>\!0,$ the matrix ${D_\a\!=\!(d_{ij}(\a))}$ defined by Eqs.\ \eqref{Qa}--\eqref{Ha1} exists and
\eqs*{
%\label{e_dist1}%-
d_{ij}(\a)=
\begin{cases}
\g\,(\a-1)                        \log_\a\!\frac{\sqrt{f_{ii}(\a)\,f_{jj}(\a)}}{f_{ij}(\a)}, &\a\ne1\vspace{0.3em}\\
\phantom{(\a-1)\,}\hspace{1.1em}\g\ln    \!\frac{\sqrt{f_{ii}( 1)\,f_{jj}( 1)}}{f_{ij}( 1)}, &\a=1
\end{cases},\quad i,j=\1n.
}
\end{prop}

\proof
Applying the matrix forest theorem~\eqref{e_mft} to $G_\a$ one obtains that $Q_\a$ exists and its entries are strictly positive, provided that $G$ is connected. Therefore $H_\a$ and $D_\a$ also exist.

Let $q_{ij}(\a)$ and $h_{ij}(\a)$ be the notation for the entries of $Q_\a$ and $H_\a$, respectively.
For any positive $\a\ne1$ and $\g$, equations \eqref{Qa}--\eqref{Da} and the matrix forest theorem~\eqref{e_mft} imply
\eqs*{
d_{ij}(\a)
&=&\tfrac{1}{2}(h_{ii}(\a)+h_{jj}(\a))-h_{ij}(\a)\\
&=&\g\,(\a-1)\left[\tfrac{1}{2}(\log_\a q_{ii}(\a)+\log_\a q_{jj}(\a))-\log_\a q_{ij}(\a)\right]\\
&=&\g\,(\a-1)\log_\a\!\frac{\sqrt{q_{ii}(\a)\,q_{jj}(\a)}}{q_{ij}(\a)}
 = \g\,(\a-1)\log_\a\!\frac{\sqrt{f_{ii}{(\a)}\,f_{jj}{(\a)}}}{f_{ij}{(\a)}}
}\nopagebreak
for every $i,j=\1n$. If $\a=1,$ then the desired expression follows similarly using~\eqref{Ha1}.
\qed
\medskip

{\noindent\bf Proof of Theorem~\ref{th_dist}.}
Proving this theorem amounts to showing that for every $i,j,k\in V(G)$:

(i) $d_{ij}(\a)=0$ if and only if $i=j$ and

(ii) $d_{ij}(\a)+d_{jk}(\a)-d_{ki}(\a)\ge0$ (triangle inequality).

\smallskip
Note that the symmetry and non-negativity of $D_\a$ (which are sometimes considered as part of the definition of distance) follow from (i) and (ii) by putting $k=j$ and $k=i$ in the triangle inequality.

Let $\a\ne1$. If $i=j,$ then by \eqref{Da}, $d_{ij}(\a)=0$. Conversely, if $d_{ij}(\a)=0,$ then by Proposition~\ref{p_d1}, $f_{ii}(\a)\,f_{jj}(\a)=(f_{ij}(\a))^2$ holds. If $i\ne j,$ then $f_{ij}(\a)<f_{jj}(\a)$, since, by definition, $\Fij(G_\a)\subseteq\Fjj(G_\a)$ and $\Fjj(G_\a)\smallsetminus\Fij(G_\a)$ contains the trivial spanning rooted forest having no edges and weight unity. Since $Q_\a$ is symmetric, $f_{ij}(\a)<f_{ii}(\a)$. Consequently, $i\ne j$\, contradicts the assumption $d_{ij}(\a)=0$, hence $i=j$.

\smallskip
To prove (ii), observe that \eqref{Ha}, \eqref{Da}, and \eqref{e_mft} for any positive $\a\ne1$ imply
\eqss{
\label{e_trine}
d_{ij}(\a)+d_{jk}(\a)-d_{ki}(\a)
&=&\tfrac{1}{2}(h_{ii}(\a)+h_{jj}(\a)+h_{jj}(\a)+h_{kk}(\a)-h_{kk}(\a)-h_{ii}(\a))\nonumber\\
&&\!\!-\;h_{ij}(\a)-h_{jk}(\a)+h_{ki}(\a)\nonumber\\
&=&h_{jj}(\a)+h_{ki}(\a)-h_{ij}(\a)-h_{jk}(\a)\nonumber\\
\label{e_tria}
&=&\g\,(\a-1)\log_\a\frac{f_{jj}(\a)\,f_{ki}(\a)}{f_{ij}(\a)\,f_{jk}(\a)}.
}
Since $Q_\a$ is symmetric and the matrix $F_\a=(f_{ij}(\a))_{n\times n}$ determines a transitional measure for $G_\a$ \cite[item~1 of Corollary~3]{Che10bottl}, we have\footnote{In this proof, we cannot formally apply Theorem~1 of \cite{Che10bottl}, since the construction of logarithmic distances in the present paper has some difference from that in~\cite{Che10bottl}.} $f_{jj}(\a)\,f_{ki}(\a)\ge f_{ij}(\a)\,f_{jk}(\a)$ (the transition inequality) and so \eqref{e_trine} implies that $d_{ij}(\a)+d_{jk}(\a)-d_{ki}(\a)\ge0$. For $\a=1$ (i) and (ii) are proved similarly.
\qed

\bigskip
Theorem~\ref{th_dist} enables us to give the following definition.

\begin{defin}
\label{def_dist}
{\rm
Suppose that $G$ is a connected weighted multigraph and $\a>0.$ A \emph{logarithmic forest distance with parameter $\a$ on $G$\/} is a function $d_\a\!:V(G)\!\times\!V(G)\to\R$ such that $d_\a(i,j)=d_{ij}(\a),$ where %$i,j\in V(G)$ and
$D_\a=(d_{ij}(\a))$ is defined by Eqs.\ \eqref{Qa}--\eqref{Ha1}.
}
\end{defin}

In Definition~\ref{def_dist}, the scaling factor $\g$ of \eqref{Ha} and \eqref{Ha1} and the transformation $G\to G_\a$ are regarded as internal parameters of logarithmic forest distances. In Section~\ref{s_bottle}, we show that all such distances are graph-geodetic. In Sections~\ref{s_asymp} and~\ref{s_weighted}, logarithmic forest distances with specific~$\g$ and $G\to G_\a$ transformations and desirable asymptotic properties are considered. Section~\ref{s_weighted} also contains natural requirements a $G\to G_\a$ transformation should satisfy.

\smallskip
\section{The logarithmic forest distances are graph-geodetic}
\label{s_bottle}

The key property of the logarithmic forest distances is that they are graph-geodetic.\footnote{This term is borrowed from~\cite{KleinZhu98}.}

\begin{defin}
\label{d_g-d}
{\rm
For a multigraph $G,$ a function $d\!:V(G)\times V(G)\to\R$ is \emph{graph-geodetic\/} whenever for all $i,j,k\in V(G),\;$ $d(i,j)+d(j,k)=d(i,k)$ holds if and only if every path in $G$ from $i$ to $k$ contains~$j$.}
\end{defin}

If $d(i,j)$ is a distance on the set of graph vertices, then the property of being graph-geodetic is a natural condition of strengthening the triangle inequality to equality. The shortest path distance clearly possesses the ``if'' (but not the ``only if'') part of the graph-geodetic property; the ``if'' part of this property for the resistance distance is established by Lemma~E in~\cite{KleinRandic93}. The ordinary distance in a Euclidean space satisfies a similar condition resulting from substituting ``line segment'' for ``path in $G$'' in Definition~\ref{d_g-d}.

\begin{thm}
\label{th_botadd}
For any connected multigraph $G$ and any {$\a>0,$} each logarithmic forest distance $d_\a(i,j)$ is graph-geodetic.
\end{thm}

Note that Theorem~\ref{th_botadd} is not tantamount to item~2 of Corollary~3 in~\cite{Che10bottl}, since the construction of logarithmic distances in the present paper differs from that in~\cite{Che10bottl}.

\bigskip
\proof
Since $F_\a=(f_{ij}(\a))_{n\times n}$ is symmetric and determines a transitional measure for $G_\a$ \cite[item~1 of Corollary~3]{Che10bottl}, we have that $f_{jj}(\a)\,f_{ki}(\a)=f_{ij}(\a)\,f_{jk}(\a)$ is true if and only if every path in $G_\a$ from $i$ to $k$ contains~$j$\/ (the graph bottleneck identity). Owing to \eqref{e_tria} and the analogous expression for $\a=1$, this equality is equivalent to ${d_\a(i,j)+d_\a(j,k)-d_\a(k,i)=0}$. On the other hand, $G_\a$ is constructed in such a way that it shares the set of paths with~$G$. Consequently, ${d_\a(i,j)+d_\a(j,k)-d_\a(k,i)=0}$ holds if and only if every path in $G$ from $i$ to $k$ contains~$j$.
\qed

\medskip\smallskip
Graph-geodetic functions have many interesting properties. One of them, as mentioned in \cite{KleinRandic93}, is a simple connection (such as that obtained in \cite{GrahamHoffmanHosoya77}) between the cofactors and the determinant of $G$'s distance matrix and those of the maximal blocks of $G$ that have no cut points. Another example is the recursive Theorem~8 in~\cite{KleinZhu98}. Clearly, for a tree, all the ${n(n-1)/2}$ values of a graph-geodetic distance are determined by the $n-1$ values corresponding to the pairs of adjacent vertices. The logarithmic forest distances, as well as their limiting cases, i.e., the shortest-path, weighted shortest-path, and resistance distances (see Sections~\ref{s_asymp} and~\ref{s_weighted}), need not be Euclidean; however, by Blumenthal's ``Square-Root'' theorem, the corresponding ``square-rooted'' distances satisfy the 3-Euclidean condition (cf.~\cite{KleinZhu98}).

It can be observed that the ``ordinary'' forest distances \cite{CheSha01} defined without the logarithmic transformation \eqref{Qa} are not generally graph-geodetic.

\section{The shortest-path and resistance distances in the framework of logarithmic forest distances}
\label{s_asymp}

Consider the family of logarithmic forest distances determined by the $G\to G_\a$ edge weight transformation
\eq{
\label{e_aunifo}%+
w_{ij}^p(\a)=\a w_{ij}^p,\;\;i,j=\1n,\;\;p=\1n_{ij}
}
(which implies $L_\a=\a L$) and the scaling factor
\eq{
\label{eq_graviy}%+
\g=\ln(e+\atn),
}
where $e$ is Euler's constant.
It turns out that the shortest-path and the resistance distances are the limiting functions of this family.

\begin{prop}
\label{p_da0}
For any connected multigraph $G$ and every $i,j\in V(G),$ $d_\a(i,j)$ with $G\!\to\! G_\a$ transformation \eqref{e_aunifo} and scaling factor~\eqref{eq_graviy} converges to the shortest path distance $d^s(i,j)$ as $\a\to0^+$.
\end{prop}

\proof
Denote by $m$ the shortest path distance $d^s(i,j)$ between $i$ and $j\ne i$.
Observe that the weight of every forest that belongs to $\Fii(G_\a)$ and has at least one edge vanishes with $\a\to0^+$, whereas $\Fii(G_\a)$ contains one forest without edges whose weight is unity.
Taking this into account and using Proposition~\ref{p_d1} and \eqref{e_fij} one obtains
\eqs*{
 \lim_{\a\to 0^+}d_\a(i,j)
=\lim_{\a\to 0^+}\left(-\log_\a\frac{\sqrt{1\cd\!1}}{\a^m\big(f_{ij}^{(m)}+o(1)\big)}\right),
}
where $o(1)\to 0$ as $\a\to 0^+$. Consequently,
$$
 \lim_{\a\to 0^+}d_\a(i,j)
=\lim_{\a\to 0^+}\big(m+\log_\a f_{ij}^{(m)}\big)=m=d^s(i,j).\eqno{\qed}
$$

\begin{prop}
\label{p_dainf}
For any connected multigraph $G$ and every $i,j\in V(G),$ $d_\a(i,j)$ with $G\!\to\! G_\a$ transformation \eqref{e_aunifo} and scaling factor~\eqref{eq_graviy} converges to the resistance distance $d^r(i,j)$ as $\a\to\infty$.
\end{prop}

\proof
Observe that for every $i,j\in V(G),\;$ $f_{ij}^{(n-1)}$ is the total weight of all spanning trees in~$G$. Denote this weight by $t$; since $G$ is connected, $t>0$. By Proposition~\ref{p_d1} one has
\eqs*{
\lim_{\a\to\infty}\!d_\a(i,j)
=\!\lim_{\a\to\infty}\!\!\left(\!\!\frac{2\a}{n}
\ln\a\,(\ln\a)^{-1}\ln\!\frac{\sqrt{\a^{n-1}\!\!\left(t+\ainv f_{ii}^{(n-2)}+o\big(\ainv\big)\!\right)
                                    \a^{n-1}\!\!\left(t+\ainv f_{jj}^{(n-2)}+o\big(\ainv\big)\!\right)\!\!}}
                                   {\a^{n-1}\!\!\left(t+\ainv f_{ij}^{(n-2)}+o\big(\ainv\big)\!\right)}
\right)\!\!,
}
where $o\big(\ainv\big)$ denotes expressions such that $\a\cd\!o\big(\ainv\big)\to0$ as $\a\to\infty$. Hence
\eqss{
\lim_{\a\to\infty}\!d_\a(i,j)
\!&=&\!\frac{2}{n}\lim_{\a\to\infty}\ln\frac{\sqrt{    \Bigl(1+\frac{f_{ii}^{(n-2)}}{\a t}\!\Bigr)^{\!\a}
                                                       \Bigl(1+\frac{f_{jj}^{(n-2)}}{\a t}\!\Bigr)^{\!\a}}}
                                                  {    \Bigl(1+\frac{f_{ij}^{(n-2)}}{\a t}\!\Bigr)^{\!\a}}
   =   \frac{2}{n}                  \ln\frac{\sqrt{\exp\Bigl(  \frac{f_{ii}^{(n-2)}}{t}\Bigr)
                                                   \exp\Bigl(  \frac{f_{jj}^{(n-2)}}{t}\Bigr)}}
                                                  {\exp\Bigl(  \frac{f_{ij}^{(n-2)}}{t}\Bigr)}\nonumber\\\nonumber
%\label{e_toinf}%-
\!&=&\!\frac{f_{ii}^{(n-2)}+f_{jj}^{(n-2)}-2f_{ij}^{(n-2)}}{nt}.
}

Consequently, by Corollary~\ref{cor_fores} of Section~\ref{sec_intro}, $\lim_{\a\to\infty}\!d_\a(i,j)=d^r(i,j)$.
\qed

\medskip
Note that for logarithmic forest distances with arbitrary positive scaling factors $\g$, ``converges'' in Propositions~\ref{p_da0} and~\ref{p_dainf} must be replaced by ``becomes proportional.''

\section{The weighted shortest path distance in the present framework}
\label{s_weighted}

In the theory of electrical networks, the edge weight $w_{ij}^p$ is interpreted as the conductance, and the Laplacian matrix $L=(\ell_{ij})$ is termed the admittance matrix. The \emph{weighted shortest path distance\/} $d^{ws}(i,j)$ is defined as follows:\footnote{This formula corrects Eq.~(6.2) in~\cite{KleinRandic93}; cf.\ the first inequality in \cite[p.~261]{DezaDeza09}.}
\eqs*{
%\label{e_dws}%-
d^{ws}(i,j)=\min_\pi\sum_{\e\in\pi}r_\e,
}
where the minimum is taken over all paths $\pi$ from $i$ to $j$ and the sum is over all edges $\e$ in $\pi$; $r_\e=1/w_\e$ is called the \emph{resistance\/} of the edge $\e$, where $w_\e$ is the weight of this edge.

It turns out that the weighted shortest path distance, as well as the ordinary shortest path distance, fits into the framework of logarithmic forest distances. To show this, it suffices to consider the $G\to G_\a$ transformation
\eq{
\label{e_aexpo}%+
w_{ij}^p(\a)=\psi_\a(r_{ij}^p),\;\;\text{where\ \ }r_{ij}^p=1/w_{ij}^p,\,\;\;i,j=\1n,\;\;p=\1n_{ij},
}
with
\eq{
\label{f_wsp}%+
\psi_\a(r)=\a^r.
}

\begin{prop}
\label{p_da00}
For any connected multigraph $G$ and every $i,j\in V(G),$ $d_\a(i,j)$ with ${G\!\to\! G_\a}$ transformation \eqref{e_aexpo}--\eqref{f_wsp}
converges to the weighted shortest path distance $d^{ws}(i,j)$  as ${\a\to0^+,}$
provided that the scaling factor $\g$ in \eqref{Ha} goes to~$1$ as $ \a\to0^+$.
\end{prop}

\proof
Let $G_\a$ be the multigraph with edge weights $\a^{r_{ij}^p}.$
Using the notation \eqref{e_fija}, for every vertices $i$ and $j\ne i$, just as in the proof of Proposition~\ref{p_da0} we derive
\eqs*{
%\label{e_affor}%-
 \lim_{\a\to 0^+}d_\a(i,j)
=\lim_{\a\to 0^+}\left(-\log_\a\frac{\sqrt{1\cd\!1}}{f_{ij}(\a)}\right)
=\lim_{\a\to 0^+}       \log_\a f_{ij}(\a).
}
For every $0<\a<1$, %it follows that %we have
\eq{
\label{e_kappa}%+
f_{ij}(\a)
=\sum_{\Ff\in\Fij(G_\a)}w(\Ff)
=\sum_{\Ff\in\Fij(G)}\prod_{\e\in E(\Ff)}\a^{r_\e}
=\sum_{\Ff\in\Fij(G)}\a^{\sum_{\e\in E(\Ff)}r_\e}
=\kappa_{ij}(\a)\ssi \a^{d^{ws}(i,j)},
}
where $1\le\kappa_{ij}(\a)\le|\Fij(G)|$.
In \eqref{e_kappa} we use the fact that for every path from $i$ to $j$, $\Fij(G)$ contains a forest sharing the set of edges with this path.
Consequently,
$$
 \lim_{\a\to 0^+}d_\a(i,j)
=\lim_{\a\to 0^+}\log_\a\big(\kappa_{ij}(\a)\ssi\a^{d^{ws}(i,j)}\big)
=d^{ws}(i,j).\eqno{\qed}
$$

\noindent
{\bf Remark.}
By definition, $G_\a$ results from $G$ by a certain transformation of edge weights. This means that $V(G_\a)=V(G)$ and for every $i,j\in V(G)$, $G$ and $G_\a$ have the same number of edges incident to both $i$ and~$j$ (this fact is used in the proof of Theorem~\ref{th_botadd}). Since the weight of every edge is positive, $\psi_\a(r)$ must be positive for every $r>0$ and every $\a$ in the definition domain. Moreover, recall that the edge weight $w_{ij}^p$ is interpreted as the conductance of the corresponding edge, $w_{ij}=\sum_{p=1}^{n_{ij}}w_{ij}^p$, and $w_{ij}=0=w_{ij}(\a)$ holds if and only if $n_{ij}=0$. Since the  absence of direct $ij$-connections, i.e., the case of $n_{ij}=0$, can also be interpreted as the zero conductance of such connections, $w_{ij}(\a)=\sum_{p=1}^{n_{ij}}w_{ij}^p(\a)$ should be \emph{small\/} whenever the conductances $w_{ij}^p$ of $ij$-edges are small (i.e., whenever their resistances are large). Formally, the continuity condition we have just described is stated as follows: $\lim_{r\to\infty}\psi_\a(r)=0$ for every $\a$ in the definition domain. Finally, a natural requirement is that $\psi_\a(r)$ should be decreasing for every~$\a$ (this monotonicity condition along with the above limiting condition implies the positivity of~$\psi_\a(r)$). Note that the %edge weight
transformation \eqref{f_wsp} satisfies these conditions if and only if $\a\in\,]0,1[$. Furthermore, the edge weight transformations we consider in this paper are such that for each $r>0,$ $\lim_{\a\to0^+}\psi_\a(r)=0$ and $\psi_\a(r)$ strictly increases in~$\a;$ except for $\psi_\a(r)=\a^r,$ these transformations are increasing functions of~$\frac{\a}{r}.$

\bigskip
Using Propositions~\ref{p_dainf} and~\ref{p_da00} one can easily define a parametric family of logarithmic forest distances converging to the weighted shortest path distance as $\a\to 0^+$ and to the resistance distance as $\a\to\infty$.
Such a family is not unique. Perhaps, the most interesting family with such asymptotic properties is the one described in Proposition~\ref{p_ws-r}.

\begin{prop}
\label{p_ws-r}
For any connected multigraph $G$ and every $i,j\in V(G),$ the logarithmic forest distance $d_\a(i,j)$ defined by\/$:$
\begin{itemize}
\Up{.4}\item Eqs.~\eqref{Qa}--\eqref{Da} with \eqref{Ha} replaced by
\Up{.4}\eq{
\label{Ha+}%+
H_\a=\g\,\a\,\ovec{\ln Q_\a},
}
\Up{2.2}\item $G\!\to\! G_\a$ transformation \eqref{e_aexpo} with
\Up{.4}\eq{
\label{f_wspr}%+
\psi_\a(r)=\dfrac{\a}{r}\,e^{-\tfrac{r}{\a}},\;\;\text{and}
}
\Up{1.8}\item any positive scaling factor $\g\!=\!\g(\a)$ such that\/ $\lim_{\a\to0^+}\g(\a)\!=\!1$ and\/ $\lim_{\a\to\infty}\g(\a)\!=\!\frac{2}{n}$
\end{itemize}
\Up{.3}converges to the weighted shortest path distance $d^{ws}(i,j)$ as $\a\to0^+$ and to the resistance distance $d^r(i,j)$ as $\a\to\infty$.
\end{prop}
%-%
%-%\begin{prop}
%-%\label{p_ws-r}
%-%For any connected multigraph $G$ and every $i,j\in V(G),$ the logarithmic forest distance $d_\a(i,j)$ defined by Eqs.~\eqref{Qa}--\eqref{Da} with \eqref{Ha} replaced by
%-%\eq{
%-%\label{Ha+}
%-%H_\a=\g\,\a\,\ovec{\ln Q_\a}
%-%}
%-%with $G\!\to\! G_\a$ transformation \eqref{e_aexpo} where
%-%\eq{
%-%\label{f_wspr}
%-%\psi_\a(r)=\dfrac{\a}{r}\,e^{-\tfrac{r}{\a}}
%-%}
%-%and with any positive scaling factor $\g\!=\!\g(\a)$ such that\/ $\lim_{\a\to0^+}\g(\a)\!=\!1$ and\/ $\lim_{\a\to\infty}\g(\a)\!=\!\frac{2}{n}$
%-%converges to the weighted shortest path distance $d^{ws}(i,j)$ as $\a\to0^+$ and to the resistance distance $d^r(i,j)$ as $\a\to\infty$.
%-%\end{prop}

Comparing \eqref{Ha+} with \eqref{Ha} shows that the family of distances introduced in Proposition~\ref{p_ws-r} is contained in the class of logarithmic forest distances \eqref{Qa}--\eqref{Ha1}. As a scaling factor in \eqref{Ha+} that meets the requirements of Proposition~\ref{p_ws-r},  one can take, for example, $\g(\a)=(\frac{2}{n}\a+\beta)/(\a+\beta)$, %or $\g(\a)=(\a+1)/(\frac{n}{2}\a+1)$.
where $\beta>0$ is a parameter.

\bigskip\proof
Let $G_\a$ be the multigraph with edge weights assigned by \eqref{e_aexpo} and \eqref{f_wspr}. Since for large $\a$, the function \eqref{f_wspr} is asymptotically equivalent to $\a/r$, using Proposition~\ref{p_dainf} we conclude that for every $i,j\in V(G),$ $\lim_{\a\to\infty}d_\a(i,j)=d^r(i,j)$.

Furthermore, for every vertices $i$ and $j\ne i$ and the distance $d_\a(i,j)$ under consideration, similarly to the proof of Proposition~\ref{p_da0} we have
\eq{
\label{e_affor1}%+
 \lim_{\a\to 0^+}d_\a(i,j)
=\lim_{\a\to 0^+}\left( \a\ln\frac{\sqrt{1\cd\!1}}{f_{ij}(\a)}\right)
=\lim_{\a\to 0^+}    (-\a\ln                      f_{ij}(\a)).
}
The definition of the graph weight and Eqs.\ \eqref{e_fija}, \eqref{e_aexpo}, and~\eqref{f_wspr} yield
\eqs*{
%\label{e_kappa1}%-
f_{ij}(\a)
=\sum_{\Ff\in\Fij(G_\a)}w(\Ff)
=\sum_{\Ff\in\Fij(G   )}\prod_{\e\in E(\Ff)}\tfrac{\a}{r_\e}e^{-r_\e/\a}
=\sum_{\Ff\in\Fij(G   )}\a^{m_\Ff}w\_\Ff e^{-d_\Ff/\a},
}
where $m\_\Ff=|E(\Ff)|$, $w\_\Ff=\prod_{\e\in E(\Ff)}w_\e,$ and $d_\Ff=\sum_{\e\in E(\Ff)}r_\e.$

Observe that if $\Ff,\Ff'\in\Fij(G)$ and                              \newline
(a) $d_{\Ff}<d_{\Ff'}$ or                                             \newline
(b) $d_{\Ff}=d_{\Ff'}$ and $m\_\Ff<m\_{\Ff'}$  or                     \newline
(c) $d_{\Ff}=d_{\Ff'}$,    $m\_\Ff=m\_{\Ff'}$, and $w\_\Ff>w\_{\Ff'},$\newline
then for each small enough $\a>0$,
$\a^{m\_\Ff}  w\_\Ff e^{-d_\Ff/\a}
>\a^{m_{\Ff'}}w\_{\Ff'} e^{-d_{\Ff'}/\a}$ holds. Consequently, there exists $\a_0>0$ such that for all $\a\in\,]0,\a_0[$
and some $\kappa_{ij}(\a)$ satisfying $1\le\kappa_{ij}(\a)\le|\Fij(G)|$,
\eq{
\label{e_kappa2}
f_{ij}(\a)
=\kappa_{ij}(\a)\,\a^{m_\Ffb}\,w\_\Ffb e^{-d^{ws}(i,j)/\a}
}
is true,
where $\Ffb$ is a forest $\Ff\in\Fij(G)$ that satisfies (a) or (b) or the nonstrict version (with $w\_\Ff\ge w\_{\Ff'}$) of (c) w.r.t.\ each
$\Ff'\in\Fij(G)$. Substituting \eqref{e_kappa2} in \eqref{e_affor1} results in
$$
  \lim_{\a\to 0^+}d_\a(i,j)
%=\lim_{\a\to 0^+}(-\a\ln\kappa_{ij}(\a)\a^{m_\Ffb}w\_\Ffb e^{-d^{ws}(i,j)/\a})
 =\lim_{\a\to 0^+}\big(-\a\big(\ln(\kappa_{ij}(\a)\,w\_\Ffb)+m_\Ffb\ln\a-d^{ws}(i,j)/\a\big)\big)
 =d^{ws}(i,j).\eqno{\qed}
$$

\section{Concluding remarks}
\label{s_conc}

Thus, the main property of the logarithmic forest distances introduced by means of Theorem~\ref{th_dist} and Proposition~\ref{p_d1} is that they are graph-geodetic: $d(i,j)+d(j,k)=d(i,k)$ if and only if every path connecting $i$ and $k$ contains~$j$ (Theorem~\ref{th_botadd}).

Three classical distances, namely, the shortest-path, the resistance, and the weighted shortest-path distances, all fit, as limiting cases, into the framework of logarithmic forest distances. The two former distances can be obtained by the use of the edge weight transformation \eqref{e_aunifo}, which generates the regularized Laplacian kernel, or, in other words, by putting $\psi_\a(r)=\a/r$ in~\eqref{e_aexpo} (Propositions~\ref{p_da0} and~\ref{p_dainf}). To obtain the latter distance, one can put\footnote{It can be shown that $\psi_\a(r)=e^{-{r}/{\a}}$ is also suitable for this purpose.} $\psi_\a(r)=\a^r$ (Proposition~\ref{p_da00}).

To define a parametric family of logarithmic forest distances whose limiting cases are the weighted shortest path distance and the resistance distance, it suffices to put $\psi_\a(r)=\tfrac{\a}{r}\,e^{-\tfrac{r}{\a}}$ in~\eqref{e_aexpo} (Proposition~\ref{p_ws-r}).

The proofs of Theorems~\ref{th_dist} and~\ref{th_botadd} are based on the fact that the matrix $F=(f_{ij})$ of spanning rooted forests determines a transitional measure~\cite{Che10bottl} on the corresponding multigraph. That is why it can be useful to study the graph-geodetic distances produced by the other transitional measures considered in~\cite{Che10bottl}.

\bigskip
We conclude with several %additional
remarks.

\vspace{-.5em}
\subsubsection*{On intercomponent distances}\vspace{-0.5em}
Throughout the paper, we assumed that $G$ is connected.
Otherwise, if $G$ has more than one component and $i$ and $j$ belong to different components, then, by the matrix forest theorem~\eqref{e_mft}, $q_{ij}=f_{ij}=0$. Consequently, if $\log_\a(\cdot)$ and $\ln(\cdot)$ are considered as functions mapping to the extended line $\R\cup\{-\infty,+\infty\},$ then \eqref{Da} leads to the distance $+\infty$ between $i$ and~$j$, which seems quite natural.

\vspace{-.5em}
\subsubsection*{On the parameter ${\bm\a}$ and the length of paths between vertices}\vspace{-0.5em}
The parameter $\a$ of logarithmic forest distances controls the relative influence of short, medium, and long paths between vertices $i$ and $j$ on the distance $d_\a(i,j)$. As $\a\to0$, only the (weighted) shortest paths matter; the long paths have the maximum effect as $\a\to\infty$.

\vspace{-.5em}
\subsubsection*{On the ``mixture'' of the shortest-path and resistance distances}\vspace{-0.5em}
The simplest way of ``generalizing'' both the (weighted) shortest-path and the resistance distances is to consider the convex combination of the form $d'_\a(i,j)=(1-\a) d^s(i,j)+{\a d^r(i,j)}$, where $\a\in[0,1]$. However, this approach seems quite poor %scanty
from both theoretical and practical points of view. First, it does not presuppose any underlying model that might provide a deeper insight by unifying the shortest-path and the resistance distances; thus, the mixture seems just ``mechanical.'' Second, consider, for example, a path on four vertices: let $V(G)=\{1,2,3,4\}$ and $E(G)=\{(1,2),(2,3),(3,4)\}$. Then $d^s(1,2)=d^s(2,3)=d^r(1,2)={d^r(2,3)=1,}$ and therefore $d'_\a(1,2)=d'_\a(2,3)$ for all $\a\in[0,1]$. On the other hand, in applications, there are models and intuitive heuristics that result in either $d(1,2)>d(2,3)$ or $d(1,2)<d(2,3)$. Indeed, suppose that the distance $d(i,j)$ should depend on the whole set of routes between $i$ and~$j$: the shorter and more numerous are the routes, the smaller must be the distance. Then the inequality $d(1,2)>d(2,3)$ is suggested by the observation that there are three routes of length~3 between vertices~2 and~3 (namely, $(2,3,2,3)$, $(2,1,2,3)$, and $(2,3,4,3)$) and only two routes of length 3 between vertices~1 and~2 ($(1,2,1,2)$ and $(1,2,3,2)$). On the other hand, if the \emph{relative\/} numbers of routes are important, then the opposite inequality $d(1,2)<d(2,3)$ can be justified by the observation that $(1,2)$ is the unique route of length~1 starting at vertex~1, whereas $(2,3)$ and $(3,2)$ are not unique routes starting at vertices 2 and 3, respectively. It is worth mentioning that the inequality $d(1,2)<d(2,3)$ holds true for the \emph{quasi-Euclidean graph distance\/}~\cite{KleinZhu98}. %,IvanciucsKlein01

The above example demonstrates that distances providing $d(1,2)=d(2,3)$ are insufficient for the numerous applications of graph theory. As regards the forest distances, the logarithmic forest distances provide $d_\a(1,2)<d_\a(2,3)$, whereas with the ``ordinary'' forest distances \cite{CheSha01}, we have $\tilde d_\a(1,2)>\tilde d_\a(2,3)$.

\vspace{-.5em}
\subsubsection*{On some physical and probabilistic interpretations of graph distances}
\vspace{-0.5em}
In the view of H.\,Chen and F.\,Zhang \cite{ChenZhang07}, ``...the shortest-path [distance] might be imagined to be more relevant when there is corpuscular communication (along edges) between two vertices, whereas the resistance distance might be imagined to be more relevant when the communication is wave-like.'' The authors do not develop this idea in depth; presumably they have in mind that a corpuscle always takes a shortest path between vertices, while a wave takes all paths simultaneously. As has been shown in this paper, the shortest-path, weighted shortest-path, and resistance distances are extreme examples of the logarithmic forest distances. The forest distance between vertices $i$ and $j$\/ is interpreted as the probability of choosing a forest partition separating $i$ and $j$ in the model of random forest partitions \cite[Proposition~5]{CheSha01}. As $\a\to0$, transformation \eqref{Ha} preserves only those partitions that connect $i$ and $j$ by a (weighted) shortest path and separate all vertices this path does not involve; thereby the (weighted) shortest path distance results in this case, as we see in Propositions~\ref{p_da0}, \ref{p_da00}, and~\ref{p_ws-r}. When $\a\to\infty$ and $\psi_\a(r)\sim\a/r$, this transformation preserves only the partitions determined by two disjoint trees, which leads to the resistance distance, as Propositions~\ref{p_dainf} and~\ref{p_ws-r} demonstrate.

%The matrices $Q_\a$ are stochastic and they have a Markov chain interpretation (given, e.\,g., in Section~4 of \cite{Che08DAM}).

\section*{Acknowledgements}

This work was partially supported by RFBR Grant 09-07-00371 and the RAS Presidium Program ``Development of Network and Logical Control.''

I am grateful to Marco Saerens for discussing with me the problem of generalizing the shortest-path and the resistance distances.


\begin{thebibliography}{10}
\providecommand{\BIBentrySTDinterwordspacing}{\spaceskip=0pt\relax}
\relax

\bibitem{BuckleyHarary90}
F.~Buckley, F.~Harary, {Distance in Graphs}, Addison-Wesley, Redwood City, CA, 1990.

%-%\bibitem{Chan69}
%-%S.P.~Chan, {Introductory Topological Analysis of Electrical Networks}, Holt, Rinehart and Winston, New York, 1969.
%-%
\bibitem{ChandraRaghavan89}
A.K.~Chandra, P.~Raghavan, W.L.~Ruzzo, R.~Smolensky, P.~Tiwari,
  The electrical resistance of a graph captures its commute and cover times,
  in: Proc. 21st Annual ACM Symp. on Theory of Computing, ACM Press, Seattle, 1989, pp.~574--586.

\bibitem{Che08DAM}
P.~Chebotarev, Spanning forests and the golden ratio, {Discrete Appl. Math.} 156 (2008) 813--821.

\bibitem{Che10bottl}
P.~Chebotarev, The graph bottleneck identity, Adv. in Appl. Math. URL: {\small\verb"http://dx.doi.org/10.1016/j.aam.2010.11.001"} (in press).
%-%\url{http://dx.doi.org/10.1016/j.aam.2010.11.001}
%-%[arXiv preprint math.CO/1003.3904, 2010. URL \verb"<http://arxiv.org/abs/1003.3904>"]. %\url{}
\BIBentrySTDinterwordspacing

\bibitem{CheAga02ap}
P.~Chebotarev, R.~Agaev, Forest matrices around the {Laplacian} matrix, {Linear Algebra Appl.} 356 (2002) 253--274.

\bibitem{CheSha95}
P.Yu.~Chebotarev, E.~Shamis,
  On the proximity measure for graph vertices provided by the inverse {Laplacian} characteristic matrix,
  in: {Abstracts of the conference ``Linear Algebra and its Applications,''} 10--12 July, 1995,
  University of Manchester, Manchester, UK, 1995, pp.~6--7. URL: {\small\verb"http://www.ma.man.ac.uk/~higham/laa95/abstracts.ps"}. %\url{}

\bibitem{CheSha97}
P.Yu.~Chebotarev, E.V.~Shamis, The matrix-forest theorem and measuring relations in small social groups, {Autom. Remote Control} 58 (1997) 1505--1514.

\bibitem{CheSha98}
P.Yu.~Chebotarev, E.V.~Shamis, On proximity measures for graph vertices, {Autom. Remote Control} 59 (1998) 1443--1459.

\bibitem{CheSha01}
P.~Chebotarev, E.~Shamis, The forest metrics for graph vertices, {Electron. Notes Discrete Math.} 11 (2002) 98--107.

\bibitem{ChenZhang07}
H.~Chen, F.~Zhang, Resistance distance and the normalized {Laplacian} spectrum, {Discrete Appl. Math.} 155 (2007) 654--661.

\bibitem{ChungZhao10}
F.~Chung, W.~Zhao, PageRank and random walks on graphs,
in: Proceedings of ``Fete of Combinatorics and Computer Science'' Conference in honor of Laci Lov\'asz, Keszthely, Hungary, August 11--15, 2008 (in press).

\bibitem{DezaDeza09}
M.M.~Deza, E.~Deza, {Encyclopedia of Distances}, Springer, Berlin--Heidelberg, 2009.

\bibitem{DezaLaurent97}
M.M. Deza, M. Laurent, {Geometry of Cuts and Metrics}, Springer, Berlin, 1997.
%Volume 15 of Algorithms and Combinatorics,

\bibitem{GobelJagers74}
F.~G\"obel, A.A. Jagers, Random walks on graphs, {Stochastic Process. Appl.} 2 (1974) 311--336.

\bibitem{GrahamHoffmanHosoya77}
R.L.~Graham, A.J.~Hoffman, H.~Hosoya, On the distance matrix of a directed graph, {J.~Graph Theory} 1 (1977) 85--88.

\bibitem{Gurvich10}
V.~Gurvich, Metric and ultrametric spaces of resistances, {Discrete Appl. Math.} 158 (2010) 1496--1505.

\bibitem{GvishianiGurvich87}
A.D.~Gvishiani,~V.A. Gurvich, Metric and ultrametric spaces of resistances, {Russian Math. Surv.} 42 (1987) 235--236.

\bibitem{KleinRandic93}
D.J.~Klein, M.~Randi\'c, Resistance distance, {J. Math. Chem.} 12 (1993) 81--95.

\bibitem{KleinZhu98}
D.J.~Klein, H.Y.~Zhu, Distances and volumina for graphs, {J. Math. Chem.} 23 (1998) 179--195.

\bibitem{MantrachYenCallut09}
A.~Mantrach, L.~Yen, J.~Callut, K.~Fran\c{c}oisse, M.~Shimbo, M.~Saerens,
The Sum-over-paths covariance kernel: A novel covariance measure between nodes of a directed graph,
  {IEEE Trans. Pattern Anal. Machine Intelligence} 32 (2010) 1112--1126.

\bibitem{Nash-Williams59}
C.St.J.A. Nash-Williams, Random walk and electric currents in networks, {Math. Proc. Cambridge Philos. Soc.} 55 (1959) 181--194.

\bibitem{SeshuReed61}
S.~Seshu, M.B.~Reed, Linear Graphs and Electrical Networks, Addison-Wesley, Reading, MA, 1961.

\bibitem{Shapiro87}
L.W.~Shapiro, An electrical lemma, {Math. Mag.} 60 (1987) 36--38.

\bibitem{ShimboItoM09}
M.~Shimbo, T.~Ito, D.~Mochihashi, Yu.~Matsumoto,
On the properties of von {Neumann} kernels for link analysis,
  {Machine Learning} 75 (2009) 37--67.

\bibitem{SmolaKondor03}
A.J.~Smola, R.I.~Kondor,
Kernels and regularization of graphs,
in: {Proc. 16th Annual Conf. on Learning Theory}, 2003, pp.~144--158.

\bibitem{YenSaerensShimbo08}
L.~Yen, M.~Saerens, A.~Mantrach, M.~Shimbo,
  A family of dissimilarity measures between nodes generalizing both the shortest-path and the commute-time distances,
  in: {14th ACM SIGKDD Intern. Conf. on Knowledge Discovery and Data Mining}, 2008, pp.~785--793.

\end{thebibliography}
\end{document}